\documentclass[a4paper]{article}
\usepackage[english]{babel}
\usepackage{amsmath,amssymb,amsfonts}
\usepackage{mathrsfs}
\usepackage{graphicx}
\usepackage{ifpdf}
\usepackage[numbers,sort&compress]{natbib}

\ifpdf
\usepackage[bookmarks=false,%
pdfstartview=FitH,linkbordercolor={0.5 1 1},%
citebordercolor={0.5 1 0.5},unicode,%
pagebackref,hyperindex]{hyperref}%
\else
\fi

\newtheorem{theorem}{Theorem}
\newcommand{\setA}{\mathscr{A}}

\begin{document}
\date{}

\title{On accuracy of approximation of the spectral
radius by the Gelfand formula\thanks{This work was supported by
the Russian Foundation for Basic Research, project no.
06-01-00256.}}

\author{Victor Kozyakin\\[5mm]
Institute for Information Transmission
Problems\\ Russian Academy of Sciences\\ Bolshoj Karetny lane 19,
Moscow 127994 GSP-4, Russia}

\maketitle

\begin{abstract}
The famous Gelfand formula $\rho(A)=
\limsup_{n\to\infty}\|A^{n}\|^{1/n}$ for the spectral radius of a matrix
is of great importance in various mathematical constructions.
Unfortunately, the range of applicability of this formula is
substantially restricted by a lack of estimates for the rate of
convergence of the quantities $\|A^{n}\|^{1/n}$ to $\rho(A)$. In the
paper this deficiency is made up to some extent. By using the Bochi
inequalities we establish explicit computable estimates for the rate of
convergence of the quantities $\|A^{n}\|^{1/n}$ to $\rho(A)$. The
obtained estimates are then extended for evaluation of the joint
spectral radius of matrix sets.

\medskip

\medskip\noindent PACS number 02.10.Ud; 02.10.Yn


\medskip

\end{abstract}

\section{Introduction}\label{S-intro}

Let $A$ be a complex $d\times d$ matrix and $\|\cdot\|$ be a
norm in ${\mathbb{C}}^{d}$. As is known, the spectral radius
$\rho(A)$ of the matrix $A$ can be expressed in terms of the
norms of its powers $\|A^{n}\|$ by the following \emph{Gelfand
formula:}
\begin{equation}\label{E-GF}
\rho(A)=
\lim_{n\to\infty}\|A^{n}\|^{1/n},
\end{equation}
which is equivalent to the equality
$$
\rho(A)= \inf_{n\ge1}\|A^{n}\|^{1/n}.
$$
Nowadays, the Gelfand formula is treated as a commonly known
fact and is mentioned  in practically all textbooks on linear
analysis without any references to the original publication,
which was apparently \cite{Gelf:MatSb41:e}.

The spectral radius of a single matrix is defined as the
maximum of modulus of its eigenvalues. For matrix sets it is
impossible to define the notion of the spectral radius in the
same manner. In this case, it is the formula (\ref{E-GF}) that
was taken in \cite{RotaStr:IM60} as the basis for the
definition of some quantity similar to the spectral radius.

Let $\setA$ be a non-empty bounded set of complex $m\times m$
matrices. As usually, for $n\ge1$ denote by $\setA^{n}$ the set
of all $n$-products of matrices from $\setA$; $\setA^{0}=I$.
Given a norm $\|\cdot\|$ in ${\mathbb{C}}^{d}$, the limit
\begin{equation}\label{E-JSRad}
\rho(\setA)=
\lim_{n\to\infty}\|\setA^{n}\|^{1/n},
\end{equation}
where
$$
\|\setA^{n}\|=\max_{A\in\setA^{n}}\|A\|=\max_{A_{i}\in\setA}
\|A_{n}\cdots
A_{2}A_{1}\|,
$$
is called the \emph{joint spectral radius} of the matrix set
$\setA$ \cite{RotaStr:IM60}. The limit in (\ref{E-JSRad})
always exists and does not depend on the norm $\|\cdot\|$.
Moreover, for any $n\ge 1$ the estimates $\rho({\setA})\le
\|\setA^{n}\|^{1/n}$ hold   \cite{RotaStr:IM60}, and therefore
the joint spectral radius can be defined also by the following
formula:
\begin{equation}\label{E-JSR}
\rho(\setA)=
\inf_{n\ge 1}\|\setA^{n}\|^{1/n}.
\end{equation}

Since for singleton matrix sets $\setA=\{A\}$ the equality
(\ref{E-JSRad}) coincides with the Gelfand formula (\ref{E-GF})
then (\ref{E-JSRad}) is sometimes called  the generalized
Gelfand formula \cite{ShihWP:LAA97}. There are also a number of
different definitions
\cite{DaubLag:LAA92,DaubLag:LAA01,ChenZhou:LAA00,Maesumi:CDC05,Prot:FU98,Prot:CDC05-1,ParJdb:LAA08}
of an analog of the spectral radius for matrix sets.

In various situations it is important to know the conditions
under which $\rho({\setA})>0$. As can be seen, for example,
from the following inequality
\begin{equation}\label{E-Bineq}
\|\setA^{d}\|\le C_{d}\,\rho(\setA)\|\setA\|^{d-1},
\end{equation}
see  \cite[Thm. A]{Bochi:LAA03}, $\rho(\setA)=0$ if and only if
$\setA^{d}=\{0\}$, that is if and only if the matrix set
$\setA$ is nilpotent.

In the case of singleton matrix sets $\setA=\{A\}$, as is shown
in a plenty of standard courses of linear analysis, the
condition $\rho(A)\neq0$ implies the inequalities
\begin{equation}\label{E-matgen}
    \gamma^{(1+\ln n)/n}\|A^{n}\|^{1/n}\le\rho(A)\le\|A^{n}\|^{1/n}
\end{equation}
with some constant $\gamma\in(0,1)$. In \cite[Lem.
2.3]{Wirth:IJRNC98} the inequalities (\ref{E-matgen}) have been
extended for the case of general matrix sets:
\begin{equation}\label{E-wirth-gen}
    \gamma^{(1+\ln
    n)/n}\|\setA^{n}\|^{1/n}\le\rho({\setA})\le\|\setA^{n}\|^{1/n}.
\end{equation}

Unfortunately, to the best of the author's knowledge, neither
exact values for $\gamma$ nor at least effectively computable
estimates for the rate of convergence of the quantities
$\|A^{n}\|^{1/n}$ and $\|\setA^{n}\|^{1/n}$ to their limits are
known. This substantially restricts the range of applicability
of the formulas (\ref{E-GF}) and (\ref{E-JSRad}). It is not
very crucial for singleton matrix sets $\setA=\{A\}$ since in
this case the value of $\rho(A)$ can be computed by other
means. However, for the case of general matrix sets the lack of
estimates for the rate of convergence of the quantities
$\|\setA^{n}\|^{1/n}$ to $\rho(\setA)$ is much more critical
since in this case, as far as is known to the author, any
alternative ways for evaluation of $\rho(\setA)$ until now are
not found.

In the paper this deficiency is made up to some extent. By
using the Bochi inequalities (\ref{E-Bineq}) we establish below
explicit computable estimates for the rate of convergence of
the quantities $\|\setA^{n}\|^{1/n}$ to $\rho(\setA)$.
Apparently, these estimates are new even for the case of matrix
families consisting of a single matrix.

The paper is organized as follows. In Introduction we have
presented a concise survey of publications related to the
problem of evaluation of the joint (generalized) spectral
radius. In Section~\ref{S-mainT} the main result of the paper,
Theorem~\ref{T-main}, is formulated. This theorem provides
explicit upper and lower bounds for the spectral radius of the
matrix set $\setA$. The proof of the main theorem is relegated
to Section~\ref{S-proofmain}, while Section~\ref{S-constC} is
devoted to evaluation of the Bochi constant $C_{d}$ playing the
key role in the main theorem.

\section{Main Theorem}\label{S-mainT}

The aim of this section is to obtain explicit estimates for the
spectral radius of a finite matrix family. The next result from
\cite[Thm. A]{Bochi:LAA03} is of principal importance in all
further considerations.

\par\addvspace{\topsep}\par\noindent\textbf{Theorem A (J. Bochi)}
\emph{Given $d\ge1$, there exists $C_{d} > 1$ such that, for
every bounded set $\setA$ of complex $d\times d$ matrices and
every norm $\|\cdot\|$ in $\mathbb{C}^{d}$,}
\begin{equation}\label{E-Bineq1}
\|\setA^{d}\|\le C_{d}\,\rho(\setA)\|\setA\|^{d-1}.
\end{equation}

In \cite{Bochi:LAA03} the value of the constant $C_{d}$ is
given only for the case $r=1$, that is when the matrix family
$\setA$ consists of a single matrix. However, intermediate
constructions from \cite{Bochi:LAA03} contain all the
information needed to find $C_{d}$. This will allow to get in
Section~\ref{S-constC} an explicit expression for $C_{d}$.

Due to the Bochi theorem, if $\rho(\setA)=0$ then
$\setA^{d}=\{0\}$, that is the matrix set $\setA$ is nilpotent.
By (\ref{E-JSR}) a converse statement is also valid:
$\setA^{d}=\{0\}$ implies $\rho(\setA)=0$. So, theoretically
verification of the condition $\rho(\setA)=0$ may be fulfilled
in a finite number of steps: it suffices only to check that all
$d$-products of matrices from $\setA$ vanish. Of course this
remark is hardly suitable in practice since even for moderate
values of $d = 3, 4$, $r= 5, 6$ the computational burden of
calculations becomes too high. Nevertheless, in what follows we
will study only the case when
$$
\rho(\setA)\neq0\quad\textrm{or, equivalently,}\quad
\setA^{d}\neq\{0\}.
$$

\begin{theorem}\label{T-main} Given $d\ge2$, for every bounded set
$\setA$ of complex $d\times d$ matrices and every norm
$\|\cdot\|$ in $\mathbb{C}^{d}$,
\begin{equation}\label{E-mainineq}
C_{d}^{-\sigma_{d}(n)/n}
\left(\frac{\|\setA\|^{d}}{\|\setA^{d}\|}\right)^{-\nu_{d}(n)/n}
\|\setA^{n}\|^{1/n}\le
\rho(\setA)\le\|\setA^{n}\|^{1/n},\qquad n=1,2,\ldots~,
\end{equation}
where
\begin{align*}
C_{d}&=\begin{cases}
2^{d}-1&\textrm{for~}r=1,\\
d^{3d/2}&\textrm{for~}r>1,
\end{cases}\\
\sigma_{d}(n)&=\begin{cases}
\frac{1}{2}\left(\frac{\ln n}{\ln 2}+1\right) \left(\frac{\ln n}{\ln
2}+2\right)&\textrm{for~}d=2,\\
\frac{(d-1)^{3}}{(d-2)^{2}}\cdot n^{\frac{\ln(d-1)}{\ln
d}}&\textrm{for~}d>2,
\end{cases}\\
\nu_{d}(n)&=\begin{cases}
\frac{\ln n}{\ln 2}+1&\textrm{for~}d=2,\\
\frac{(d-1)^{2}}{d-2}\cdot n^{\frac{\ln(d-1)}{\ln d}}&\textrm{for~}d>2.
\end{cases}
\end{align*}
\end{theorem}

The proof of Theorem~\ref{T-main} is relegated to
Section~\ref{S-proofmain}. Clearly, the statement of
Theorem~\ref{T-main} holds also for real matrix sets.

Note that the estimates (\ref{E-mainineq}) are weaker than the
estimates (\ref{E-wirth-gen}). It is not clear now whether it
is caused by the techniques of proof of the estimates
(\ref{E-mainineq}) or by the fact that the obtained constants
$C_{d}$, $\sigma_{d}(n)$ and $\nu_{d}(n)$ are universal, that
is depend neither on a matrix set nor on the choice of the norm
$\|\cdot\|$.

Note also that the value of the constant $C_{d}$ rapidly
increases in $d$. That is why the estimates (\ref{E-mainineq})
are hardly useful in applications and sooner are of theoretical
interest. Moreover, the estimates (\ref{E-mainineq}) are
essentially finite-dimensional and scarcely can be extended for
linear operators in infinite-dimensional spaces.

Remark, at last, that for irreducible matrix sets $\setA$
containing more that one matrix there are valid \cite[Lem.
2.3]{Wirth:IJRNC98} the following, stronger than
(\ref{E-wirth-gen}) or (\ref{E-mainineq}), estimates:
$$
    \gamma^{1/n}\|\setA^{n}\|^{1/n}\le\rho(\setA)\le\|\setA^{n}\|^{1/n},
$$
where the constant $\gamma$ can be effectively computed
\cite{Koz:ArXiv08-2}.

\section{Proof of Theorem~\ref{T-main}}\label{S-proofmain}

The inequality $\rho(\setA)\le\|\setA^{n}\|^{1/n}$ in
(\ref{E-mainineq}) follows from (\ref{E-JSR}). For $r=1$ the
value of the constant $C_{d}$ is found in \cite{Bochi:LAA03};
for $r>1$ this constant will be evaluated in
Section~\ref{S-constC}.

Let us deduce some corollaries from the Bochi theorem. Firstly
note that for any natural numbers $p$ and $q$ the following
inequalities hold
\begin{equation}\label{E-auxineq}
\|\setA^{p+q}\|\le\|\setA^{p}\|\cdot\|\setA^{q}\|,
\end{equation}
from which
\begin{equation}\label{E-auxineq1}
\|\setA^{p}\|\le\|\setA\|^{p},\quad
\rho(\setA^{p})=\rho^{p}(\setA),\qquad p=1,2,\ldots~.
\end{equation}
Then from (\ref{E-Bineq}) we immediately get:
$$
\|\setA^{d^{k}}\|\le C_{d}\left(\rho(\setA)\right)^{d^{k-1}}
\|\setA^{d^{k-1}}\|^{d-1},\qquad k=1,2,\ldots~.
$$
If we denote
$$
\omega_{n}(\setA)=
\frac{\|\setA^{{n}}\|}{\left(\rho(\setA)\right)^{{n}}},\qquad
n=1,2,\ldots~,
$$
then the latter inequalities can be rewritten in the form:
$$
\omega_{d^{k}}(\setA)\le
C_{d}\left(\omega_{d^{k-1}}(\setA)\right)^{d-1},\qquad k=1,2,\ldots~.
$$
Therefore, for any integer $k=1,2,\ldots$
\begin{align*}
\omega_{d^{k}}(\setA)&\le
C_{d}\left(\omega_{d^{k-1}}(\setA)\right)^{d-1},\\
\left(\omega_{d^{k-1}}(\setA)\right)^{d-1}&\le
C_{d}^{d-1}\left(\omega_{d^{k-2}}(\setA)\right)^{(d-1)^{2}},\\
\left(\omega_{d^{k-2}}(\setA)\right)^{(d-1)^{2}}&\le
C_{d}^{(d-1)^{2}}\left(\omega_{d^{k-3}}(\setA)\right)^{(d-1)^{3}},\\
&\ldots\\
\left(\omega_{d}(\setA)\right)^{(d-1)^{k-1}}&\le
C_{d}^{(d-1)^{k-1}}\left(\omega_{1}(\setA)\right)^{(d-1)^{k}}.
\end{align*}
By multiplying the obtained inequalities we get:
\begin{equation}\label{E-omegadn}
\omega_{d^{k}}(\setA)\le C_{d}^{\sum_{i=0}^{k-1}(d-1)^{i}}
\left(\omega_{1}(\setA)\right)^{(d-1)^{k}},\qquad k=1,2,\ldots~.
\end{equation}
Now, note that by the Bochi inequality (\ref{E-Bineq1})
$$
\frac{1}{\rho(\setA)}\le C_{d}\frac{\|\setA\|^{d-1}}{\|\setA^{d}\|}.
$$
Hence
$$
1\le\omega_{1}(\setA)= \frac{\|\setA\|}{\rho(\setA)}\le
C_{d}\frac{\|\setA\|^{d}}{\|\setA^{d}\|}.
$$
This allows to derive from (\ref{E-omegadn}) the estimate for
$\omega_{d^{k}}(\setA)$ which does not contain in the
right-hand part the unknown value $\rho(\setA)$:
\begin{equation}\label{E-omegadn-fin}
\omega_{d^{k}}(\setA)\le C_{d}^{\sum_{i=0}^{k}(d-1)^{i}}
\left(\frac{\|\setA\|^{d}}{\|\setA^{d}\|}\right)^{(d-1)^{k}},\qquad
k=0,1,\ldots~.
\end{equation}

Now, let $n$ be an arbitrary natural number. Then there is a
natural $k$ such that
$$
d^{k}\le n <d^{{k+1}},
$$
and consequently for $n$ it is valid the representation
$$
n=n_{k}d^{k}+n_{k-1}d^{k-1}+\dots+n_{0},
$$
where
\begin{equation}\label{E-kn}
1\le n_{k}\le d-1,\qquad 0\le n_{i}\le d-1,\quad i=1,2,\ldots,k-1.
\end{equation}

Since by (\ref{E-auxineq}) and (\ref{E-auxineq1})
$$
\omega_{p+q}(\setA)\le\omega_{p}(\setA)\cdot\omega_{q}(\setA)
$$
for any natural numbers $p$ and $q$, then
$$
\omega_{n}(\setA)\le\left(\omega_{d^{k}}(\setA)\right)^{n_{k}}\cdot
\left(\omega_{d^{k-1}}(\setA)\right)^{n_{k-1}}\cdots
\left(\omega_{1}(\setA)\right)^{n_{0}}.
$$
By (\ref{E-omegadn-fin}) from here it follows
\begin{equation}\label{E-omega-k}
\omega_{n}(\setA)\le C_{d}^{\sigma_{d}(n)}
\left(\frac{\|\setA\|^{d}}{\|\setA^{d}\|}\right)^{\nu_{d}(n)},
\end{equation}
where
\begin{equation}\label{E-signu}
\sigma_{d}(n)=\sum_{j=0}^{k}n_{j}\sum_{i=0}^{j}(d-1)^{i},\qquad
\nu_{d}(n)=\sum_{j=0}^{k}n_{j}(d-1)^{j}.
\end{equation}

Note that, by definition of the value $\omega_{n}(\setA)$,
(\ref{E-omega-k}) is equivalent to
$$
\|\setA^{n}\|\le C_{d}^{\sigma_{d}(n)}
\left(\frac{\|\setA\|^{d}}{\|\setA^{d}\|}\right)^{\nu_{d}(n)}
\left(\rho(\setA)\right)^{n},
$$
and therefore to the inequality
$$
C_{d}^{-\sigma_{d}(n)/n}
\left(\frac{\|\setA\|^{d}}{\|\setA^{d}\|}\right)^{-\nu_{d}(n)/n}
\|\setA^{n}\|^{1/n}\le \rho(\setA).
$$

Since this last inequality coincides with the left-hand part of
(\ref{E-mainineq}) then to complete the proof of the theorem it
remains only to get the estimates for $\sigma_{d}(n)$ and
$\nu_{d}(n)$. By (\ref{E-kn}) and (\ref{E-signu})
\begin{multline}\label{E-sigma-k}
\sigma_{d}(n)=\sum_{j=0}^{k}n_{j}\sum_{i=0}^{j}(d-1)^{i}\le
(d-1)\sum_{j=0}^{k}\sum_{i=0}^{j}(d-1)^{i}=\\
(d-1)\sum_{j=0}^{k}(k+1-j)(d-1)^{j},
\end{multline}
\begin{equation}
\label{E-nu-k}
\nu_{d}(n)=\sum_{j=0}^{k}n_{j}(d-1)^{j}\le
(d-1)\sum_{j=0}^{k}(d-1)^{j}.
\end{equation}

By definition of the number $k$ we have $k\le\frac{\ln n}{\ln
d}$. Then for $d=2$ from (\ref{E-sigma-k}), (\ref{E-nu-k}) it
follows:
\begin{alignat*}{2}
\sigma_{2}(n)&\le\frac{(k+1)(k+2)}{2}~&\le~&
\frac{1}{2}\left(\frac{\ln
n}{\ln 2}+1\right) \left(\frac{\ln n}{\ln 2}+2\right),\\
\nu_{2}(n)&\le k+1~&\le~& \frac{\ln n}{\ln 2}+1.
\end{alignat*}

Represent (\ref{E-sigma-k}), (\ref{E-nu-k}) for $d>2$ in the
form
\begin{alignat}{2}\label{E-sigma-k3}
\sigma_{d}(n)&=\sum_{j=0}^{k}n_{j}\sum_{i=0}^{j}(d-1)^{i}~
&\le~&
(d-1)^{k+1}\sum_{j=0}^{k}\frac{j+1}{(d-1)^{j}},\\
\label{E-nu-k3} \nu_{d}(n)&=\sum_{j=0}^{k}n_{j}(d-1)^{j}~
&\le~& (d-1)^{k+1}\sum_{j=0}^{k}\frac{1}{(d-1)^{j}},
\end{alignat}
and use the equalities
$$
\sum_{j=0}^{\infty}x^{j}=\frac{1}{1-x},\quad
\sum_{j=0}^{\infty}(j+1)x^{j}=\frac{1}{(1-x)^{2}},\qquad |x|<1.
$$
By setting here $x=\frac{1}{d-1}$, from (\ref{E-sigma-k3}),
(\ref{E-nu-k3}) we obtain:
\begin{alignat*}{2}
\sigma_{d}(n)&\le\frac{(d-1)^{k+3}}{(d-2)^{2}}~ &\le~
&\frac{(d-1)^{3}}{(d-2)^{2}}\cdot n^{\frac{\ln(d-1)}{\ln d}},\\
\nu_{d}(n)&\le\frac{(d-1)^{k+2}}{d-2}~ &\le~
&\frac{(d-1)^{2}}{d-2}\cdot n^{\frac{\ln(d-1)}{\ln d}}.
\end{alignat*}

The theorem is proved.

\section{Evaluation of $\boldsymbol{C_{d}}$}\label{S-constC}

In \cite{Bochi:LAA03} existence of the constant $C_{d}$ is
established in Theorem A, proof of which is based on Lemmas 2
and 3 cited below.

\par\addvspace{\topsep}\par\noindent\textbf{Lemma 2 (J. Bochi)}
\emph{Let $\|\cdot\|_{e}$ be the Euclidian norm in
$\mathbb{C}^{d}$. There exists $C_{0}=C_{0}(d)$ such that
$$
\|S\setA^{d}S^{-1}\|_{e}\le
C_{0} \|\setA\|_{e}\|S\setA S^{-1}\|_{e}^{d-1}
$$
for every non-empty bounded set $\setA$ of $d\times d$ matrices
and every matrix $S\in GL(d)$.}\par\addvspace{\topsep}\par

Actually, in \cite{Bochi:LAA03} under the proof of Lemma 2 it
is obtained first that for every diagonal matrix $S\in GL(d)$
the following inequality holds
$$
\|S\setA^{d}S^{-1}\|_{0}\le
d^{d-1} \|\setA\|_{0}\|S\setA S^{-1}\|_{0}^{d-1}.
$$
with the matrix norm $\|A\|_{0}=\max|a_{ij}|$.

As is known \cite[Ch. 5]{HJ:e}, the following relations between
the norm $\|\cdot\|_{0}$ and the Euclidean norm $\|\cdot\|_{e}$
hold:
$$
\|A\|_{0}\le \|A\|_{e}\le d \|A\|_{0},
$$
from which the chain of inequalities follows:
\begin{multline*}
d^{-1}\|S\setA^{d}S^{-1}\|_{e}\le \|S\setA^{d}S^{-1}\|_{0}\le
d^{d-1} \|\setA\|_{0}\|S\setA S^{-1}\|_{0}^{d-1}\le\\
\|S\setA^{d}S^{-1}\|_{0}\le d^{d-1} \|\setA\|_{e}\|S\setA
S^{-1}\|_{e}^{d-1},
\end{multline*}
that is
$$
\|S\setA^{d}S^{-1}\|_{e}\le
d\cdot
d^{d-1} \|\setA\|_{e}\|S\setA S^{-1}\|_{e}^{d-1}.
$$
The last inequality, as shown in \cite{Bochi:LAA03} under the
proof of Lemma 2, can be easily extended to the general case
$S\in GL(d)$. Therefore $C_{0}=d^{d}$.

Now, let us move to consideration of Lemma 3 from
\cite{Bochi:LAA03}.

\par\addvspace{\topsep}\par\noindent\textbf{Lemma 3 (J. Bochi)}
\emph{There exists $C=C(d)$ such that, for every two norms
$\|\cdot\|_{1}$ and $\|\cdot\|_{2}$ in $\mathbb{C}^{d}$ there
is a matrix $S\in GL(d)$ such that}
\begin{enumerate}
\item $C^{-1}\|v\|_{1}\le \|Sv\|_{2}\le\|v\|_{1}$ \emph{for
    all} $v\in\mathbb{C}^{d}$;
\item $C^{-1}\|A\|_{1}\le \|SAS^{-1}\|_{2}\le C\|A\|_{1}$
    \emph{for all $d\times d$ matrices} $A$.
\end{enumerate}\par\addvspace{\topsep}\par

Here the second part is an immediate consequence of the first
one. To evaluate the constant $C$ in the first part, first
notice that whenever Lemma 3 is applied in \cite{Bochi:LAA03},
one of the two norms $\|\cdot\|_{1}$ or $\|\cdot\|_{2}$ is the
Euclidian norm.

So, let us evaluate the constant $C$ under the assumption that
the norm $\|\cdot\|_{1}$ is arbitrary while the norm
$\|\cdot\|_{2}$ is Euclidean. This can be done by using a
matrix-theoretic version of complex John's ellipsoid theorem
\cite{AndoShih:SIAM:98}. Certainly J.~Bochi was not aware of
this technique when he wrote his paper. To be more specific,
let us reproduce the argumentation from \cite{Shih:TM99}.

Given a norm $\|\cdot\|_{1}$ in $\mathbb{C}^{d}$, it can be
represented in the form
$$
\|v\|_{1}^{2}= \sup_{\lambda\in\Lambda}\langle H_{\lambda}v,v\rangle,
\quad v\in\mathbb{C}^{d},
$$
where $\langle\cdot,\cdot\rangle$ is the Euclidean scalar
product in $\mathbb{C}^{d}$ and $\{H_{\lambda},
\lambda\in\Lambda\}$ is a family of semidefinite matrices. But
according to \cite[Thm.~2.1]{AndoShih:SIAM:98} for any family
of semidefinite matrices $\{H_{\lambda}, \lambda\in\Lambda\}$
there is a positive definite matrix $H$ such that
$$
\langle Hv,v\rangle \le
\sup_{\lambda\in\Lambda}\langle H_{\lambda}v,v\rangle \le
d\langle Hv,v\rangle, \quad v\in\mathbb{C}^{d}.
$$
Therefore
$$
\langle Hv,v\rangle \le
\|v\|_{1}^{2} \le
d\langle Hv,v\rangle, \quad v\in\mathbb{C}^{d}.
$$
Since the matrix $H$ may be thought of as symmetric then, by
setting $S=H^{1/2}$,
$\|\cdot\|_{2}=\sqrt{\langle\cdot,\cdot\rangle}$ and
$\|Sv\|_{2}^{2}=\langle Sv,Sv\rangle\equiv\langle
H^{1/2}v,H^{1/2}v\rangle\equiv \langle Hv,v\rangle$, we obtain
$$
d^{-1}\|v\|_{1}^{2}\le \|Sv\|_{2}^{2}\le\|v\|_{1}^{2},
$$
and the conclusion of Lemma 3 is valid with the constant $C =
d^{1/2}$.

Now, to evaluate the value of the constant $C_{d}$ in Theorem A
it suffices to note that due to \cite{Bochi:LAA03}
$C_{d}=C^{d}C_{0}$ where $C_{0}$ and $C$ are the constants from
Lemmas 2 and 3, respectively. Hence, $C_{d}=d^{3d/2}$.

\section*{Acknowledgments}
I am grateful to Jairo Bochi for many valuable comments and
especially for the idea to use John's ellipsoid theorem for
estimating the constant $C$ in Lemma 3 which have allowed to
improve dramatically the estimate of $C_{d}$.

 \ifdefined\inputencoding\inputencoding{cp1251}\fi \newcommand{\nosort}[1]{}
  \newcommand{\bbljan}[0]{January} \newcommand{\bblfeb}[0]{February}
  \newcommand{\bblmar}[0]{March} \newcommand{\bblapr}[0]{April}
  \newcommand{\bblmay}[0]{May} \newcommand{\bbljun}[0]{June}
  \newcommand{\bbljul}[0]{July} \newcommand{\bblaug}[0]{August}
  \newcommand{\bblsep}[0]{September} \newcommand{\bbloct}[0]{October}
  \newcommand{\bblnov}[0]{November} \newcommand{\bbldec}[0]{December}

\end{document}